\documentclass[a4paper,12pt]{amsproc}
\usepackage{amsmath,amsfonts,amssymb,amsthm,anysize}
\usepackage{enumerate}
\usepackage{epsfig}
\usepackage{supertabular}
\usepackage{graphicx}
\usepackage{color}

\newcommand{\tr}{\text{Tr}}

\newcommand\Z{{\mathbb Z}}
\newcommand\Q{{\mathbb Q}}
\newcommand\R{{\mathbb R}}
\newcommand\T{{\mathrm{Tr}}}
\newcommand\C{{\mathbb C}}
\newcommand\F{{\mathbb F}}

\newcommand\Gal{{\mathrm{Gal}}}

\newcommand\Tr{{\mathrm{Tr}}}

\newcommand\Cay{{\mathrm{Cay}}}
\newcommand\ord{{\mathrm{ord}}}

\newcommand\PG{\mathsf{PG}}
\newcommand\QQ{{\mathsf{Q}}}
\newcommand\HH{{\mathsf{H}}}

\renewcommand\mod{{\mathrm{mod\, \, }}}
\newcommand\pro{{\bf Proof: }}

\newcommand\cQ{{\mathcal Q}}

\newcommand\cP{{\mathcal P}}

\theoremstyle{plain}
\newtheorem{theorem}{Theorem}[section]
\newtheorem{problem}[theorem]{Problem}
\newtheorem{lemma}[theorem]{Lemma}
\newtheorem{corollary}[theorem]{Corollary}
\newtheorem{fact}[theorem]{Fact}

\newtheorem{proposition}[theorem]{Proposition}

\numberwithin{equation}{section}

\theoremstyle{remark}
\newtheorem{remark}[theorem]{Remark}

\usepackage{url, hyperref}
\hypersetup{citecolor=blue, linkcolor=blue, colorlinks=true}

\def\id{{\mathrm{id}}}
\def\pe{{\mathfrak{p}}}
\def\Pe{{\mathfrak{P}}}
\renewcommand\le{\leqslant}
\renewcommand\ge{\geqslant}

\begin{document}

\title[Cyclic arcs and strongly regular graphs]{Cyclic arcs of Singer type and strongly regular Cayley graphs over finite fields}

\author{Koji Momihara$^*$}
\address{ %
	Division of Natural Science\\
	Faculty of Advanced Science and Technology\\
	Kumamoto University\\
	2-40-1 Kurokami, Kumamoto 860-8555, Japan}
\email{momihara@educ.kumamoto-u.ac.jp}
\thanks{$^*$Research supported by 
	JSPS under Grant-in-Aid for Scientific Research (C) 20K03719.}

\author{Qing Xiang$^\dagger$}
\address{ %
	Department of Mathematics and National Center for Applied Mathematics Shenzhen\\
	Southern University of Science and Technology\\
	Shenzhen 518055, China
}
\email{xiangq@sustech.edu.cn}
\thanks{$^\dagger$Research partially supported by the National Natural Science Foundation of China grant 12071206, and the Sino-German Mobility Programme M-0157}


\begin{abstract}
In \cite{M18}, the first author gave a construction  of strongly regular Cayley graphs on the additive group of finite fields by using three-valued Gauss periods. 
In particular, together with the result in \cite{BLMX}, it was shown 
that there exists a strongly regular Cayley  graph 
with negative Latin square type parameters 
$(q^6,r(q^3+1),-q^3+r^2+3r,r^2+r)$, where  
$r=M(q^2-1)/2$, in the following cases: 
(i) $M=1$ and $q\equiv 3\,(\mod{4})$; 
(ii) $M=3$ and $q\equiv 7\,(\mod{24})$; and (iii) $M=7$ and  $q\equiv 11,51\,(\mod{56})$. The existence of strongly regular Cayley graphs with the above parameters for odd $M>7$  was left open. In this paper, we prove that 
if there is an $h$, $1\le h\le M-1$, such that $M\,|\,(h^2+h+1)$ and the order of $2$ in $(\Z/M\Z)^\times$ is odd, 
then there exist infinitely many primes $q$ such that strongly regular Cayley graphs 
with the aforementioned parameters exist. 
\end{abstract}


\maketitle

\section{Introduction}\label{sec:intro}

A {\it strongly regular graph} srg$(v,k,\lambda,\mu)$ is a simple and undirected graph, neither complete nor edgeless, that has the following properties:

\begin{itemize}

\item [(1)] It is a regular graph of order $v$ and valency $k$.

\item [(2)] For each pair of adjacent vertices $x,y$, there are exactly $\lambda$ vertices adjacent to both $x$ and $y$.

\item [(3)] For each pair of nonadjacent vertices $x,y$, there are exactly $\mu$ vertices adjacent to both $x$ and $y$.
\end{itemize}

For example, the pentagon is an srg$(5,2,0,1)$ and the Petersen graph is an srg$(10,3,0,1)$. The following spectral characterization of srgs is well known. Let $\Gamma$ be a simple connected $k$-regular graph that is neither complete nor edgeless. Then $\Gamma$  is strongly regular if and only if its adjacency matrix has exactly two distinct eigenvalues different from $k$. 

An srg$(u^2,r(u-\epsilon),\epsilon u+r^2-3\epsilon r,r^2-\epsilon r)$ is said to be of 
{\it Latin square type or negative Latin square type} according as $\epsilon=1$ or $-1$. Typical examples of strongly regular graphs of Latin square type or negative  Latin square type come from nonsingular 
quadrics in the projective space $\PG(n,q)$ of odd dimension $n$ over the finite field $\F_q$ of order $q$~\cite{Ma}. It seems that more examples of 
strongly regular  graphs of Latin square type are known than  those of negative Latin square type~\cite{DX04}.

Let $G$ be an additively written abelian group, and let $D$ be an inverse-closed subset of $G$ such that  $0_G\not\in D$. The {\it Cayley graph on $G$ with connection set $D$}, denoted by ${\rm Cay}(G,D)$, is the graph with the elements of $G$ as vertices; two vertices are adjacent if and only if their difference belongs to $D$.  It is known that the eigenvalues of $\Cay(G,D)$ are the character values of $D$. More precisely, for a (complex) character $\psi$ of $G$, define 
\[
\psi(D)=\sum_{x\in D}\psi(x). 
\] 
Then all the eigenvalues of $\Cay(G,D)$ are given by $\psi(D)$, $\psi\in {\hat G}$, where ${\hat G}$ is the (complex) character group of $G$. 
  
During the past two decades, strongly regular Cayley graphs on the additive groups of finite fields have been extensively studied due to their close connections with certain substructures in finite geometries. If the connection set $D$ of $\Cay(\F_{q^n},D)$ is invariant under multiplication by elements of the multiplicative group $\F_{q}^\ast$, then the subset ${\mathcal D}=\{x\F_{q}^\ast: x\in D\}$ of $\F_{q^n}^\ast/\F_{q}^\ast$ can be viewed as a set of points of $\PG(n-1,q)$; in this case, $\Cay(\F_{q^n},D)$ is strongly regular if and only if ${\mathcal D}$ is a projective two-intersection set, or a two-character set, in $\PG(n-1,q)$. Further links between strongly regular Cayley graphs and substructures in finite classical polar spaces such as $m$-ovoids and $i$-tight sets are known~\cite{BKLP,BLP,FWX}. In this paper we will only consider strongly regular Cayley graphs related to $m$-ovoids. Let ${\mathcal S}$ be a finite classical polar space and ${\mathcal M}$ be a set of points in ${\mathcal S}$. We say that ${\mathcal M}$ is an {\it $m$-ovoid} if every generator of ${\mathcal S}$ meets ${\mathcal M}$ in exactly $m$ points. There are many known constructions of $m$-ovoids~\cite{BLP,CCEM,CP1,CP3,D,E,FMX4,FWX,K}. 

An elliptic quadric $\QQ^-(5,q)$ of $\PG(5,q)$ is trivially a $(q+1)$-ovoid, which gives rise to a negative Latin square type  srg$(q^6,r(q^3+1),-q^3+r^2+3 r,r^2+ r)$ with $r=q^2-1$. 
On the other hand, it is known that nontrivial $m$-ovoids in $\QQ^-(5,q)$ exist whenever $q$ is an odd prime power and $m=(q+1)/2$. Since the dual of $\QQ^-(5,q)$ as a generalized quadrangle is $\HH(3,q^2)$, such a $\frac{(q+1)}{2}$-ovoid in $\QQ^-(5,q)$ can be interpreted as a hemisystem in $\HH(3,q^2)$, which consists of $\frac{(q+1)(q^3+1)}{2}$ lines containing exactly half of the lines 
through every point of $\HH(3,q^2)$. Hemisystems have been studied extensively in \cite{BGR,BLMX,CP2,CP,KNS,S}. They give rise to negative Latin square type srg$(q^6,r(q^3+1),-q^3+r^2+3 r,r^2+ r)$ with $r=(q^2-1)/2$. In \cite{BLMX}, the authors constructed a new $\frac{(q+1)}{2}$-ovoid  for $q\equiv 3\,(\mod{4})$ in  $\QQ^-(5,q)$ by choosing $2(q+1)$ cosets from the $4(q^2+q+1)$th cyclotomic classes of $\F_{q^6}$.  
The first author of this paper could generalize the construction of strongly regular Cayley graphs in \cite{BLMX} by using three-valued Gauss periods~\cite{M18}. In particular, the following theorem was proved. 

\begin{theorem} {\em (\cite{M18})}
Let $q$ be an odd prime power. 
There exists an srg$(q^6,r(q^3+1),-q^3+r^2+3 r,r^2+ r)$ with $r=(q^2-1)M/2$, in the following cases:  
\begin{itemize}
\item[(1)] $M=3$ and $q\equiv 7\,(\mod{24})$; 
\item[(2)] $M=7$ and $q\equiv 11,51\,(\mod{56})$. 
\end{itemize}
\end{theorem}

The three-valued Gauss periods used in \cite{M18} all arise from cyclic arcs of Singer type in $\PG(2,q)$. An $M$-{\it arc} in $\PG(2,q)$ is a set of $M$ points no three of which are collinear. An arc is said to be {\it cyclic} if it is a point orbit under the action of a cyclic collineation group $H$ of $\PG(2,q)$ on the point set of $\PG(2,q)$. If $H$ is a subgroup of a Singer group of $\PG(2,q)$ the arc is said to be of {\it Singer type}. It turns out that besides the cyclic arcs of Singer type of size $3$ and $7$ used in \cite{M18}, there are a lot more ``small" cyclic arcs of Singer type known \cite{M95, SZ96} in $\PG(2,q)$. A natural question arises: Can these cyclic $M$-arcs of Singer type with $M>7$ give rise to strongly regular Cayley graphs? In this paper, we investigate this problem and come up with sufficient conditions for these cyclic $M$-arcs of Singer type to give rise to srg$(q^6,r(q^3+1),-q^3+r^2+3 r,r^2+ r)$ with $r=(q^2-1)M/2$ and $M>7$. For the detailed statement of our results, see Theorem~\ref{thm:const_sum}.

We further investigate to what extent our construction in Theorem~\ref{thm:const_sum} works. Let $M,h$ be positive integers such that $M$ is odd,  $1\le h\le M-1$ and  $M\,|\,(h^2+h+1)$, and let ${\mathcal P}_{M,h}$ be the set of primes $p\equiv h\,(\mod{M})$. We use 
$\Psi_{M,h}$ to denote the set of primes $p\equiv h\,(\mod{M})$ with  
$p\equiv 3\,(\mod{4})$ such that there exists 
a strongly regular Cayley graph on $(\F_{p^6},+)$  
with parameters $(p^{6},r(p^3+1),-p^3+r^2+3r,r^2+r)$, where $r=(p^2-1)M/2$. 
Furthermore, for the quadratic character $\eta$ of $\F_{p^3}$ and a fixed primitive element $\omega\in \F_{p^3}$, define 
\[
\Psi_{M,h,\alpha,\beta}=
\{p\in {\mathcal P}_{M,h}\,:\,
\eta(1+\omega^\frac{(p^3-1)i}{M})=\alpha, 1\le i\le M-1, \eta(2)=-\alpha,
\eta(-1)=\beta \mbox{ in } \F_{p^3}\},
\]
where $\alpha,\beta\in \{1,-1\}$. In this paper, we will prove that $\Psi_{M,h}$ contains almost all primes in $\bigcup_{\alpha\in \{1,-1\}}\Psi_{M,h,\alpha,-1}$. 
In particular, the following theorem is proved. 

\begin{theorem}\label{main:ext_thm1}
If  the order of $2$ in $(\Z/M\Z)^\times$ is odd,  then $\Psi_{M,h}$ contains infinitely many primes. 
\end{theorem}

The values of $M<200$ satisfying the condition in Theorem~\ref{main:ext_thm1} 
are 
\[
7,31,49,73,79,103,127,151,199.
\] 
Thus, using the $M$'s listed above and Theorem~\ref{main:ext_thm1} we obtain many new infinite families of strongly regular Cayley graphs of negative Latin square type.
%

The organization of this paper is as follows. In Section~\ref{sec:pre}, we review a 
construction of strongly regular Cayley graphs based on three-valued Gauss periods given in \cite{M18}. Furthermore, we state Chebotar\"{e}v's density theorem, which will be used to prove Theorem~\ref{main:ext_thm1}. In Section~\ref{section:infinite}, we give a lower bound on primes $p$ such that the Gauss periods in $\F_{q^3}$ of characteristic $p$ are guaranteed to take exactly three values forming an arithmetic progression. This result is needed for applying our construction of strongly regular Cayley graphs in Theorem~\ref{thm:const_sum}. In Section~\ref{sec:cond}, we show that 
$\Psi_{M,h}$ contains almost all primes in $\bigcup_{\alpha\in\{1,-1\}}\Psi_{M,h,\alpha,-1}$. Furthermore, we study conditions under which $\Psi_{M,h,\alpha,\beta}$ is an infinite set. In particular, we determine 
$\Psi_{M,h,\alpha,\beta}$ in the cases where $M=3,7$, or $21$. In Section~\ref{sec:inf}, we study the structure of the Galois group $\Gal(E_M/\Q)$, where $E_M$ is obtained  by adjoining $\sqrt{1+\zeta_M^i}$, $0\le i\le M-1$, and $\zeta_4$ to  $\Q$. Finally we give a proof of Theorem~\ref{main:ext_thm1}. In Section~\ref{sec:conclu}, we propose a  few open problems for future work.


\section{Preliminaries}\label{sec:pre}
\subsection{Cyclotomic classes, Gauss periods, and  Gauss sums}
Let $p$ be a prime and let $\zeta_p=\exp(2\pi i/p)\in \C$ be a primitive $p^{\rm th}$ root of unity. For positive integers  $f$ and $n$, let $q=p^f$ and let $\F_{q^n}$ denote the finite field of order $q^n$. 
Define 
$$\psi_{\F_{q^n}}:\F_{q^n}\to \C^\ast, \;\psi_{\F_{q^n}}(x)=\zeta_p^{\Tr_{q^n/p}(x)},\;\forall x\in \F_{q^n}$$ 
where $\T_{q^n/p}$ is the trace map from $\F_{q^n}$ to $\F_p$. It can be easily shown that $\psi_{\F_{q^n}}$ is a nontrivial character of the additive group of $\F_{q^n}$. All characters of $(\F_{q^n},+)$ are given by $\psi_a$, $a\in \F_{q^n}$, where $\psi_{a}(x)=\psi_{\F_{q^n}}(ax)$, $\forall x\in \F_{q^n}$. 

Let $\omega$ be a fixed  primitive  element of $\F_{q^n}$ and $N$ be a positive integer dividing $q^n-1$. The cosets  $C_i^{(N,q^n)}=\omega^i\langle \omega^N\rangle$, $0\le i\le N-1$, of $\langle \omega^N\rangle$ in $\F_{q^n}^*$ are called the {\it $N^{\rm th}$ cyclotomic classes} of $\F_{q^n}$. The numbers $(i,j)_N:=|(C_i^{(N,q^n)}+1)\cap C_j^{(N,q^n)}|$, $0\le i,j\le N-1$, are called {\it cyclotomic numbers}.

The {\it $N^{\rm th}$ Gauss periods} of $\F_{q^n}$ are defined as the character values of cyclotomic classes: 
\[
\psi_{\F_{q^n}}(C_i^{(N,q^n)})=\sum_{x\in C_i^{(N,q^n)}}\psi_{\F_{q^n}}(x), \, \, \, 0\le i\le N-1. 
\] 
On the other hand, for a multiplicative character $\chi$ of  $\F_{q^n}$, the {\it Gauss sum} is defined by
\[
G_{q^n}(\chi)={\sum_{x\in \F_{q^n}^\ast}}\chi(x)\psi_{\F_{q^n}}(x). 
\]  
By the orthogonality relations of characters, the $N^{\rm th}$ Gauss periods of 
$\F_{q^n}$ can be expressed as a linear combination of Gauss sums:
\begin{equation}\label{eq:ortho}
\psi_{\F_{q^n}}(C_i^{(N,q^n)})=\frac{1}{N}\sum_{j=0}^{N-1}G_{q^n}(\chi_N^{j})\chi_N^{-j}(\omega^i), \; \, 0\le i\le N-1,
\end{equation}
where $\chi_N$ is a fixed multiplicative character of order $N$ of $\F_{q^n}$.  

The following is a well-known result on Eisenstein sums. 

\begin{theorem}\label{thm:Yama}{\em (\cite[Theorem~12.1.1]{BEW97})}
Let $\chi$ be a nontrivial multiplicative character of $\F_{q^n}$ whose restriction to $\F_q$ is trivial. Let $L$ be a system of coset representatives of $\F_{q}^\ast$ in $\F_{q^n}^\ast$ such that $\Tr_{q^n/q}$ maps $L$ onto $\{0,1\}\subseteq \F_p$.  Define  
$L_0=\{x\in L:\,\Tr_{q^n/q}(x)=0\}, \, {\rm and}\; L_1=L\setminus L_0.  
$ 
Then
\[
\sum_{x\in L_0}\chi(x)=\frac{1}{q}G_{q^n}(\chi). 
\]
\end{theorem}
Let $S=\{i\,(\mod{\frac{q^n-1}{q-1}}): \omega^i\in L_0\}$. Then $S$ is a Singer difference set in $\Z_{(q^n-1)/(q-1)}$, and all the hyperplanes of $\PG(n-1,q)$ can be obtained by cyclically shifting $S$. Let 
$\chi_N$ be a multiplicative character of order $N$ of $\F_{q^n}$. 
In the rest of this subsection we always assume that $N\,|\,(q^n-1)/(q-1)$.  Under this assumption, the restriction  of $\chi_N$ to $\F_q$ is trivial; by \eqref{eq:ortho} and Theorem~\ref{thm:Yama}, we have 

\begin{align}\label{eq:Gaussperio}
\psi_{\F_{q^n}}(C_i^{(N,q^n)})=&\, -\frac{1}{N}+\frac{1}{N}\sum_{j=1}^{N-1}G_{q^n}(\chi_N^{j})\chi_N^{-j}(\omega^i)\nonumber\\
=&\, -\frac{q^n-1}{N(q-1)}+\frac{q}{N}\sum_{j=0}^{N-1}\sum_{\ell\in S}\chi_N^j(\omega^{\ell-i}). 
\end{align}

Let $S_N$ denote the multiset $\rho(S)$, where $\rho$ is the natural epimorphism $\rho:\Z_{(q^n-1)/(q-1)}\rightarrow \Z_N$. We identify $S_N$ with the group ring element $S_N=\sum_{i=0}^{N-1}c_i [i]\in \Z[\Z_N]$, where 
$$c_i=|S\cap \{jN+i: j=0,1,\ldots,\frac{(q^n-1)}{N(q-1)}-1\}|,\; 0\le i\le N-1.$$ 
These numbers are the intersection sizes of $C_0^{(N,q^n)}/\F_q^*$ with the hyperplanes of $\PG(n-1,q)$. Define
\begin{equation}\label{eq:defFN}
F_N=\{c_i:0\le i\le N-1\}. 
\end{equation}
That is, $F_N$ consists of distinct intersection sizes of $C_0^{(N,q^n)}/\F_q^*$ with the hyperplanes of $\PG(n-1,q)$. Also let 
\begin{equation}\label{eq:defIIII}
I_\beta=\{i\in \Z_N:c_i =\beta\}, \, \beta \in F_N. 
\end{equation}
Then,  we have 
\begin{equation}\label{eq:sNequ}
S_N=\sum_{\beta\in F_N}\beta I_\beta \in \Z[\Z_N].
\end{equation} 
By \eqref{eq:Gaussperio}, the Gauss period $\psi_{\F_{q^n}}(C_i^{(N,q^n)})$ with $i\in I_\beta$ is related to hyperplane intersection size by 
\begin{equation}\label{eq:fiform}
\psi_{\F_{q^n}}(C_i^{(N,q^n)})= -\frac{q^n-1}{N(q-1)}+q\beta. 
\end{equation}

In \cite{FMX2}, the authors studied the question of when the Gauss periods $\psi_{\F_{q^n}}(C_i^{(N,q^n)})$, $i=0,1,\ldots,N-1$, with $N\,|\,(q^n-1)/(q-1)$ take exactly three values forming an arithmetic progression in relation to the existence of circulant weighing matrices and $3$-class translation association schemes. One important special case where the Gauss periods $\psi_{\F_{q^n}}(C_i^{(N,q^n)})$, $i=0,1,\ldots,N-1$, with $N\,|\,(q^n-1)/(q-1)$ take exactly three values forming an arithmetic progression is when $n=3$ and $C_0^{(N,q^3)}/\F_q^*$ is a cyclic arc of Singer type in $\PG(2,q)$. Both Maruta~\cite{M95} and Sz\"onyi \cite{SZ96} studied 
the question of when $C_0^{(N,q^3)}/\F_q^*$ is a cyclic arc of Singer type in $\PG(2,q)$.

Let $n=3$, $N|(q^2+q+1)$, and $M=(q^2+q+1)/N$. Assume that $C_0^{(N,q^3)}/\F_q^*$ is an $M$-arc in $\PG(2,q)$. Then each line of $\PG(2,q)$ meets $C_0^{(N,q^3)}/\F_q^*$ in $0,1$, or $2$ points. So we have $\beta_1=2$, $\beta_2=1$, $\beta_3=0$ using the notation introduced above. Hence by \eqref{eq:fiform}, the Gauss periods $\psi_{\F_{q^n}}(C_i^{(N,q^n)})$ take exactly three values $\alpha_1,\alpha_2,\alpha_3$, with $(\alpha_1,\alpha_2,\alpha_3)=(-M+2q,-M+q,-M)$. To simplify notation we write $I_i:=I_{\beta_i}$. 
Then, by Lemma~2.5 in \cite{FMX2}, we have 
\begin{equation}\label{eq:sizeI}
|I_1|=\frac{M-1}{2}, \, |I_2|=q-M+2,\, |I_3|=\frac{q^2+q+1}{M}-q+\frac{M-3}{2}. 
\end{equation}
Furthermore, by \eqref{eq:sNequ}, we have  
\begin{equation}\label{eq:quotients}
S_N=2I_1+I_2\in  \Z[\Z_{N}]. 
\end{equation}
\subsection{A construction of srgs based on three-valued Gauss periods}\label{section:construction}
Let $q\equiv 3\,(\mod{4})$ be a prime power and $\omega$ be a fixed primitive element of $\F_{q^3}$. 
Let $N$ be an odd positive integer dividing $q^2+q+1$, and let  $C_i^{(N,q^3)}=\omega^i\langle \omega^N\rangle$, $i=0,1,\ldots,N-1$. 
We suppose that the Gauss periods $\psi_{\F_{q^3}}(C_i^{(N,q^3)})$, $i=0,1,\ldots,N-1$, take exactly three values $\alpha_1=-M+2q,\alpha_2=-M+q,\alpha_3=-M$. Recall that 
\[
I_j=\{i\,(\mod{N})\,|\,{\psi_{\F_{q^3}}}(C_i^{(N,q^3)})=\alpha_j\}, \, \, \, j=1,2,3. 
\]
Let $T_1,T_2$ be a partition of  $I_2$, and 
let  
$T_i'\equiv 4^{-1}T_i\,(\mod{N})$ 
for 
$i=1,2$. Define 
\begin{equation}\label{eq:defX1}
X=2T_1' \cup (2T_2'+N)\,(\mod{2N})
\end{equation}
and 
\begin{align}
Y_X=&\, \{Ni+4j \, (\mod{4N}): (i,j)\in (\{0,3\}\times T_1')\cup   (\{1,2\}\times T_2')\}  \nonumber\\
&\, \, \, \, \cup \{Ni+4j\, (\mod{4N}): i=0,1,2,3,\, 
j\in 4^{-1}I_1\,(\mod{N})\}. 
\label{eq:defI}
\end{align} 

Let $\gamma$ be a fixed primitive element of $\F_{q^{6}}$ such that 
$\gamma^{q^3+1}=\omega$.
Define 
\begin{equation}\label{eq:dddd}
D_X=\bigcup_{i \in Y_X}C_i^{(4N,q^{6})},
\end{equation}
where $C_i^{(4N,q^{6})}=\gamma^i \langle \gamma^{4N}\rangle $, $i=0,1,\ldots,4N-1$.  It is clear that $|D_X|=(q^6-1)(2|I_1|+|I_2|)/2N=(q^2-1)(q^3+1)M/2$. 

\begin{proposition}{\em (\cite[Proposition~4.2]{M18})}\label{mainconstruction1}
For $a\in \Z_{4N}$, define $b\equiv 4^{-1}a\,(\mod{N})$ and $c\equiv 2b\,(\mod{2N})$. Let $\eta$ be the quadratic character of $\F_{q^3}$. If  $X$ defined in \eqref{eq:defX1} satisfies the condition:
\begin{align}\label{eq:3-chara}
&\, 2\psi_{\F_{q^3}}(\omega^c \bigcup_{\ell\in X}C_{\ell}^{(2N,q^3)})-
\psi_{\F_{q^3}}(\omega^c \bigcup_{\ell\in 2^{-1}I_2}C_{\ell}^{(N,q^3)})\\
=&\, \left\{
\begin{array}{ll}
\pm G_{q^3}(\eta), & \mbox{ if $c\in 2^{-1}I_2\,(\mod{N})$,}\\
0, & \mbox{ otherwise, }
 \end{array}
\right. \nonumber
\end{align}
then $\Cay(\F_{q^{6}},D_X)$ is a strongly regular graph with negative Latin square type parameters  
$(q^{6},r(q^3+1),-q^3+r^2+3r,r^2+r)$ with $r=(q^2-1)M/2$.

\end{proposition}

\subsection{A partition of a conic and its reduction modulo $2N$}\label{section:conic}
%
Viewing $\F_{q^3}$ as a $3$-dimensional vector space over $\F_q$, we will use $\F_{q^3}$ as the underlying vector space of $\PG(2,q)$. The points of $\PG(2,q)$ are $\langle \omega^i\rangle:=\omega^i \F_{q}^\ast$,
$0\le i\le q^2+q$. 
Define a quadratic form $Q: \F_{q^3}\rightarrow \F_q$ by setting $Q(x)=
\tr_{q^3/q}(x^2)$. It is straightforward to check that $Q$ is nonsingular. 
Therefore, $Q$ defines a conic $\cQ$ in $\PG(2,q)$, which contains $q+1$ points. Consequently each line of $\PG(2,q)$ meets $\cQ$ in $0$, $1$ or $2$ points. 
Consider the following subset of $\Z_{q^2+q+1}$: 
\[
W_\cQ=\{i\, (\mod{q^2+q+1}):Q(\omega^i)=0\}=\{d_0,d_1,\ldots, d_{q}\},
\]
where the elements are numbered in any unspecified order. 
Then, the conic $\cQ$ is equal to $\{\langle \omega^{d_i}\rangle:\,0\le i\le q\}$.  
Furthermore, $W_\cQ\equiv 2^{-1}S\, (\mod{q^2+q+1})$, where
$S=\{i\,(\mod{q^2+q+1}):\Tr_{q^3/q}(\omega^i)=0\}$ is a Singer difference set in $\Z_{q^2+q+1}$. Hence, 
\begin{equation}\label{eq:SingerW}
W_\cQ W_\cQ^{(-1)}=\Z_{q^2+q+1}+q[0]\in  \Z[\Z_{q^2+q+1}], 
\end{equation}
where $W_\cQ^{(-1)}=\{-x\, (\mod{q^2+q+1}):x\in W_{\cQ}\}$.

We consider a partition of  $W_\cQ$. 
For $d_0\in W_\cQ$, we define
\begin{equation}\label{eqn_defX00}
{\mathcal X}_{\mathcal Q}:=\{\omega^{d_i}\Tr_{q^3/q}(\omega^{d_0+d_i}):\,1\le i\le q\}\cup\{2  \omega^{d_0}\}
\end{equation}
and 
\[
X_{\mathcal Q}:=\{\log_{\omega}(x)\,(\mod{2(q^2+q+1)}):\, x\in {\mathcal X}_{\mathcal Q}\}\subseteq \Z_{2(q^2+q+1)}. 
\]
It is clear that $X_{\mathcal Q}\pmod{q^2+q+1}=W_\cQ$. 
The subset $X_{\mathcal Q}\subseteq \Z_{2(q^2+q+1)}$ can be written as 
\[
X_{\mathcal Q}\equiv 2E_1\cup (2E_2+(q^2+q+1))\pmod{2(q^2+q+1)}
\]
for some $E_1,E_2\subseteq \Z_{q^2+q+1}$ with $|E_1|+|E_2|=q+1$. That is, we are partitioning $X_{\mathcal Q}$ into its {\it even} and {\it odd} parts. 
It follows that  
$$W_\cQ\equiv 2(E_1\cup E_2)\,(\mod{q^2+q+1})$$ and $$S\equiv 4(E_1\cup E_2)\,(\mod{q^2+q+1}),$$ i.e., the partition of $X_{\mathcal Q}$ into its even and odd parts induces partitions of $W_\cQ$ and $S$, respectively. 
We give the following important properties of $X_{\mathcal Q}$. 
\begin{lemma} {\em (\cite[Lemma~3.4]{FMX1})} \label{lem:inva}
If we use any other $d_i$ in place of $d_0$ in the definition of ${\mathcal X}_{\mathcal Q}$, then the resulting set $X_{\mathcal Q}'$ has the property that 
$X_{\mathcal Q}'\equiv X_{\mathcal Q}$ or  $X_{\mathcal Q}+(q^2+q+1)\,(\mod{2(q^2+q+1)})$. 
\end{lemma}

\begin{proposition} {\rm(\cite[Theorem 3.7, Remark 3.8]{FMX1}, \cite[Theorem~3.4]{BLMX})} \label{thm:main2}
With notation as above, 
\begin{align*}
&\, 2\psi_{\F_{q^3}}(\omega^c \bigcup_{i\in X_{\mathcal Q}}C_i^{(2(q^2+q+1),q^3)})-\psi_{\F_{q^3}}(\omega^c \bigcup_{i\in W_{\mathcal Q}}C_i^{(q^2+q+1,q^3)})\\
=&\, \left\{
\begin{array}{ll}
\pm G_{q^3}(\eta), & \mbox{if $c\,(\mod{q^2+q+1}) \in W_\cQ$, }\\
0, & \mbox{otherwise}. 
 \end{array}
\right.
\end{align*}
where $\eta$ is the quadratic character of $\F_{q^3}$. 
\end{proposition}

\begin{remark}
Assume that $N=q^2+q+1$ (hence $M=1$). Then, $\psi_{\F_{q^3}}(C_i^{(N,q^3)})$, $i=0,1,\ldots,N-1$, take exactly two values $-M+q$ and $-M$. This can be 
viewed as the situation where $I_1=\emptyset$ and $I_2=S$ in the setting of Subsection~\ref{section:construction}.  Hence, by Propositions~\ref{mainconstruction1} and \ref{thm:main2}, $\Cay(\F_{q^{6}},D_{X_{\mathcal Q}})$ is a strongly regular graph with negative Latin square type parameters  
$(q^{6},r(q^3+1),-q^3+r^2+3r,r^2+r)$ with $r=(q^2-1)/2$. 
\end{remark}


We next consider the reduction of $X_{\mathcal Q}$ modulo $2N$. Let $N|(q^2+q+1)$ and $M=(q^2+q+1)/N$. Assume that $C_0^{(N,q^3)}/\F_q^*$ is an $M$-arc in $\PG(2,q)$. Then, the Gauss periods $\psi_{\F_{q^3}}(C_i^{(N,q^3)})$, $i=0,1,\ldots,N-1$, take exactly three 
values $-M+2q,-M+q,-M$. Since $X_{\mathcal Q}\equiv W_\cQ\equiv 2^{-1}S\, (\mod{q^2+q+1})$, the reduction of $X_{\mathcal Q}$ modulo $N$ (as a multiset) is 
$2^{-1}S_N=2^{-1}(I_1\cup I_1\cup I_2)$ as seen in \eqref{eq:quotients}. 
Define 
\begin{equation}\label{eq:defXi}
X_i=[x\,(\mod{2N}):x \in X_{\mathcal Q},x\,(\mod{N})\in 2^{-1}I_i],  \, \, i=1,2, 
\end{equation}
where we use $[\ldots]$ to indicate that the $X_i$ are multisets. Then, as multisets, $X_1\; (\mod{N})=2^{-1}(I_1 \cup I_1)$ and $X_2\; (\mod{N}) = 2^{-1}I_2$. 
We say that a multiset defined over a group $G$ is {\it purely a subset} of $G$ if each element in $G$ appears in the multiset with multiplicity at most one.  Clearly, $X_2$ defined above is purely a subset of $\Z_{2N}$, but $X_1$ may not be purely a subset of $\Z_{2N}$. 
\begin{proposition}\label{prop:AB}{\em (\cite[Proposition~5.5]{M18})}
If $X_1$ is purely a subset of $\Z_{2N}$, 
it holds that 
\begin{align}
2\psi_{\F_{q^3}}(\omega^c \bigcup_{\ell\in X_2}C_{\ell}^{(2N,q^3)})-
\psi_{\F_{q^3}}(\omega^c \bigcup_{\ell\in 2^{-1}I_2}C_{\ell}^{(N,q^3)})
=\left\{
\begin{array}{ll}
\pm G_{q^3}(\eta), & \mbox{ if $c\in 2^{-1}I_2\,(\mod{N})$,}\\
0, & \mbox{ otherwise. }
 \end{array}
\right. \nonumber
\end{align}
\end{proposition}
In summary, we have the following theorem. 

\begin{theorem}\label{thm:const_sum}
Let $q\equiv 3\,(\mod{4})$ be a prime power, $N$ be a positive integer dividing $q^2+q+1$, and $M=(q^2+q+1)/N$. Assume that  
\begin{itemize}
\item[(1)] $\psi_{\F_{q^3}}(C_i^{(N,q^3)})$, $i=0,1,\ldots,N-1$, take exactly three values $-M,-M+q,-M+2q$; and  
\item[(2)] $X_1$ is purely a 
subset of $\Z_{2N}$.
\end{itemize} 
Then $\Cay(\F_{q^6},D_{X_2})$ is a strongly regular graph with parameters $(q^{6},r(q^3+1),-q^3+r^2+3r,r^2+r)$, where  $r=(q^2-1)M/2$. Here, $X_i$, $i=1,2$, are defined in \eqref{eq:defXi} and $D_{X}$ is defined in \eqref{eq:dddd}. 
\end{theorem}
\pro 
Apply Propositions~\ref{mainconstruction1} and  \ref{prop:AB} with $X= X_2$. 
\qed

\subsection{Chebotar\"{e}v's density theorem}
Let $M$ be an odd positive integer. 
Let $h$ be a positive integer such that $M\,|\,(h^2+h+1)$, 
and   ${\mathcal P}_{M,h}$ be the set of primes such that $p\equiv h\,(\mod{M})$. For 
$\alpha,\beta\in \{1,-1\}$, 

define   
\begin{equation}\label{def:psi}
\Psi_{M,h,\alpha,\beta}=
\{p\in {\mathcal P}_{M,h}\,:\,
\eta(1+\omega^\frac{(p^3-1)i}{M})=\alpha, 1\le i\le M-1, \eta(2)=-\alpha,
\eta(-1)=\beta \mbox{ in } \F_{p^3}\},  
\end{equation}
where $\eta$ is the quadratic character of $\F_{p^3}$. 
In this paper, we will consider the question of whether there are infinitely many primes in $\Psi_{M,h,\alpha,\beta}$. In particular, we will show that $\Psi_{M,h}$ contains almost all primes in $ \bigcup_{\alpha\in\{1,-1\}}\Psi_{M,h,\alpha,-1}$, 
where $\Psi_{M,h}$ is defined in Section 1. 
To study this problem, we will use the Chebotar\"{e}v density theorem, which is a generalization of the well-known Dirichlet theorem on primes in arithmetic progressions. 

Let $F$ be a finite  Galois extension of an algebraic number field $E$. 
Let ${\mathcal O}_F$ and ${\mathcal O}_E$ be the rings of integers in $F$ and $E$, respectively. 
Let $\pe$ be a prime ideal in $E$ unramified in $F$ and $\Pe$ be a prime ideal in $F$ lying over $\pe$. Then, there is a unique monomorphism $h$
from $\Gal(({\mathcal O}_F/\Pe)/({\mathcal O}_E/\pe))$ to $\Gal(F/E)$ such that $h(\sigma)(\Pe)=\Pe$ for any $\sigma \in \Gal(({\mathcal O}_F/\Pe)/({\mathcal O}_E/\pe))$ and the map from ${\mathcal O}_F/\Pe$ to itself induced by $h(\sigma):{\mathcal O}_F\to {\mathcal O}_F$ coincides with  $\sigma$. In particular, the image $\sigma_\Pe$ of the Frobenius automorphism $x\mapsto  x^{|{\mathcal O}_E/\pe|}$ by $h$ is called the {\it Frobenius substitution} with respect to $\Pe$ in $F/E$. The Frobenius substitution 
depends on the choice of $\Pe$ lying over $\pe$ only up to conjugation, i.e., 
$\sigma_{\tau\Pe}=\tau \sigma_{\Pe}\tau^{-1}$ for any $\tau\in \Gal(F/E)$. 

Let  $P(E)$ be the set of prime ideals in $E$, and let $P'(E)$ be the set of prime ideals in $E$  unramified in $F$. For any $\sigma\in \Gal(F/E)$, we use $C_\sigma$ to denote the conjugacy class of $\sigma \in \Gal(F/E)$, i.e., 
$C_\sigma=\{\tau^{-1}\sigma\tau\,:\,\tau\in \Gal(F/E)\}$; and let $S_\sigma$ be the set of prime ideals of $E$ defined by  
\[
S_\sigma=\{\Pe \cap E\in P'(E)\,:\, \Pe \mbox{ is a prime ideal in $F$ such that $\sigma_{\Pe}\in C_\sigma$}\}. 
\] 
For any subset $S$ of $P(E)$, we define the {\it natural density} of $S$ to be
\[
\lim_{x\to \infty}\frac{|\{\pe\in S:N(\pe)\le x\}|}{|\{\pe\in P(E):N(\pe)\le x \}|}
\]
if this limit exists, where $N(\pe)$ is the absolute norm of $\pe$. If the natural density of $S$ exists, then it is actually equal to the {\it Dirichlet density} of $S$ defined by 
\[
\lim_{s\to 1^+}\Big(\sum_{\pe \in S}\frac{1}{N(\pe)^s}\Big)\Big/ \Big(\sum_{\pe\in P(E)}\frac{1}{N(\pe)^s}\Big). 
\]
The following is known as Chebotar\"{e}v's density theorem~\cite{SL96}. 

\begin{theorem}\label{thm:chev}
The density of $S_\sigma$ is equal to 
$\frac{|C_\sigma|}{|G|}$; in particular, if $|C_\sigma|\ne 0$, then there are infinitely many prime ideals in $S_\sigma$. Here the claim is valid with either notion of density. 
\end{theorem}
We denote both the Dirichlet density and the natural density of $S_\sigma$ by $\delta(S_\sigma)$. 
We will use Theorem~\ref{thm:chev} in Section~\ref{sec:inf}. 


\section{Cyclic arcs of Singer type in $\PG(2,q)$}\label{section:infinite}


In this section, we consider the question of when the first assumption in Theorem~\ref{thm:const_sum} is satisfied. That is, we consider the question of when the Gauss periods $\psi_{\F_{q^3}}(C_i^{(N,q^3)})$, $i=0,1,\ldots,N-1$, take exactly three 
values $-M+2q,-M+q,-M$. Note that the Gauss periods $\psi_{\F_{q^3}}(C_i^{(N,q^3)})$, $i=0,1,\ldots,N-1$, take exactly three 
values $-M+2q,-M+q,-M$ if and only if $C_0^{(N,q^3)}/\F_q^*$ is an $M$-arc in $\PG(2,q)$. Furthermore, it is clear that $C_0^{(N,q^3)}/\F_q^*$ is an $M$-arc in $\PG(2,q)$ if and only if $\omega^{j_1N},\omega^{j_2N},\omega^{j_3N}$ are linearly independent over $\F_q$ for all distinct $0\le j_1,j_2,j_3\le M-1$. 
Observe that 
\begin{align*}
&\mbox{ $a_1\omega^{j_1N}+a_2\omega^{j_2N}+a_3\omega^{j_3N}\not=0$ for all 
$(a_1,a_2,a_3)\in \F_q^3\setminus \{(0,0,0)\}$}\\
&\mbox{\hspace{6cm} and for all distinct $0\le j_1,j_2,j_3\le M-1$} \\
\Longleftrightarrow \hspace{0.1cm}& \mbox{ $1+b_1\omega^{i_1N}\not =b_2\omega^{i_2N}$ for all
$(b_1,b_2)\in \F_q^2\setminus \{(0,0)\}$}\\
&\mbox{\hspace{6cm} and for all distinct $1\le i_1,i_2\le M-1$} \\
\Longleftrightarrow \hspace{0.1cm}& \mbox{ $(1+C_0^{(N,q^3)})\cap  C_0^{(N,q^3)}=\F_q^{\ast}\setminus \{1\}$, i.e., $(0,0)_N=q-2$.} 
\end{align*}
Thus, there is a close relationship  between cyclic arcs of Singer type in $\PG(2,q)$ and cyclotomic numbers.  
The following is our main theorem in this section.

\begin{theorem}\label{thm:main_threevalued}
Let  $M,h$ be positive integers such that $1\le h\le M-1$ and $M\,|\,(h^2+h+1)$. 
Let $q$ be a power of a prime $p$ such that $q\equiv h\,(\mod{M})$. 
Let $N=\frac{q^2+q+1}{M}$ and $\omega$ be a primitive element of $\F_{q^3}$. Suppose that 
\[
p> \Big(
\frac{{18}M}{\phi(M)}
\Big)^{\phi(M)/2\ord_{M}(p)}, 
\]
where $\phi$ is the Euler totient function. 
Then $C_0^{(N,q^3)}/\F_q^*$ is an $M$-arc in $\PG(2,q)$. 
\end{theorem}

\begin{remark}
The theorem above was essentially proved by Maruta~\cite{M95} and Sz\"onyi \cite{SZ96}, where they did not give an explicit lower bound on $p$.  We will give a proof by using a couple of recent results in the study of cyclotomic numbers~\cite{DLS19}. 
\end{remark}

\begin{proposition}
{\em (\cite[Theorem~4.1]{DLS19})}\label{prop:cyclo1}
Let $q$ be a power of a prime $p$ and $M,e$ be positive integers such that  $q=eM+1$.  Let $\omega$ be a primitive element of $\F_{q}$, $\zeta_{M}=e^{2\pi i/M}$, and let $f(x)=\sum_{i=0}^{M-1}a_ix^i\in \Z[x]$. Suppose that 
\[
p> \Big(
\frac{M}{\phi(M)}
\sum_{i=0}^{M-1}a_i^2\Big)^{\phi(M)/2\ord_{M}(p)}. 
\]  
Then $f(\omega^e)=0$ in $\F_{q}$ if and only if 
$f(\zeta_{M})=0$ in $\C$. 
\end{proposition}
\begin{proposition}
{\em (\cite[Proposition~5.8]{DLS19})}\label{prop:cyclo2}
Let $M\ge 3$ be a positive integer. Let $s_i,t_i$, $i=1,2,3$, be integers modulo $M$ such that the $s_i$'s are pairwise distinct, the $t_i$'s are pairwise distinct, and    
the $(s_i-t_i)$'s modulo $M$ are pairwise distinct. 
Then, 
\[
\zeta_{M}^{s_1+t_2}+\zeta_{M}^{s_2+t_3}+\zeta_{M}^{s_3+t_1}
-\zeta_{M}^{s_1+t_3}-\zeta_{M}^{s_2+t_1}-\zeta_{M}^{s_3+t_2}\not=0. 
\]
\end{proposition}

\vspace{0.3cm}

{\bf Proof of Theorem~\ref{thm:main_threevalued}:}
The set  $C_0^{(N,q^3)}/\F_q^*$ of points in $\PG(2,q)$ is an $M$-arc if and only if 
\[
M_{j_1,j_2,j_3}=
\begin{pmatrix}
\omega^{j_1N} & \omega^{j_2N}  & \omega^{j_3N} \\
\omega^{j_1Nq}  & \omega^{j_2Nq}  & \omega^{j_3Nq} \\
\omega^{j_1Nq^2}  & \omega^{j_2Nq^2}  &\omega^{j_3Nq^2} \\
\end{pmatrix}
\]
is nonsingular for all distinct $j_1,j_2,j_3\in \{0,1,\ldots ,M-1\}$. 

Define 
$g_{i,j}^{(h)}(x)=x^j+x^{i(h+1)}+x^{i+j(h+1)}-x^i-x^{j(h+1)}-x^{j+i(h+1)}\in \Z[x]$, where 
the exponents $j,i(h+1),i+j(h+1),i,j(h+1),j+i(h+1)$  are reduced modulo $M$. 
Let $A_1=\{j,i(h+1),i+j(h+1)\}$ and $A_2=\{i,j(h+1),j+i(h+1)\}$, where each element is reduced modulo $M$. Assume that $i,j\not\equiv0\,(\mod{M})$ and $i\not\equiv j\,(\mod{M})$. Then, $A_1\cap A_2=\emptyset$ since $\gcd{(M,h)}=
\gcd{(M,h+1)}=1$. Furthermore, it is easily checked that if two elements in $A_1$ (resp. $A_2$) are equal, 
so are all three elements in $A_1$ (resp. $A_2$). 
Hence, by writing $g_{i,j}^{(h)}(x)=\sum_{i=0}^{M-1}a_ix^{i}$, $a_i\in \Z$, $\forall i$, we have 
$\sum_{i=0}^{M-1}a_i^2=6,12$ or $18$. 

Noting that $\omega^{N(q-1)(q+1)}=\omega^{N(q-1)(h+1)}$, we have 
\[
\det (M_{j_1,j_2,j_3})=\omega^{j_1N(q^2+q-1)+Nj_2+Nj_3}g_{j_2-j_1,j_3-j_1}^{(h)}(\omega^{N(q-1)}). 
\]
Hence $M_{j_1,j_2,j_3}$ is nonsingular if and only if  $g_{j_2-j_1,j_3-j_1}^{(h)}(\omega^{N(q-1)})\not=0$. Moreover, by applying Proposition~\ref{prop:cyclo1} with $\sum_{i=0}^{M-1}a_i^2=18$, 
$g_{j_2-j_1,j_3-j_1}^{(h)}(\omega^{N(q-1)})\not=0$ if and only if 
$g_{j_2-j_1,j_3-j_1}^{(h)}(\zeta_M)\not=0$. 

We now set 
\[
(s_1,s_2,s_3)=(0,(j_2-j_1)h,(j_2-j_1)(1+h)) 
\]
and 
\[
(t_1,t_2,t_3)=(j_3(1+h)-j_1-j_2h,j_3-j_1,j_2-j_1).  
\]
Then the $s_i$ are pairwise distinct, and  the $t_i$ are pairwise distinct. Furthermore, $s_1-t_1=-j_3(1+h)+j_1+j_2h$, $s_2-t_2=(j_2-j_1)h-(j_3-j_1)$ and $s_3-t_3=(j_2-j_1)h$, which 
are all distinct. Hence, by Proposition~\ref{prop:cyclo2}, we have 
$g_{j_2-j_1,j_3-j_1}^{(h)}(\zeta_M)\not=0$. This completes the proof of the theorem. \qed


\section{Conditions for $X_1$ to be purely a subset of $\Z_{2N}$}\label{sec:cond}

	In this section, we consider the question of when the other assumption in Theorem~\ref{thm:const_sum} is satisfied; that is, we consider when the multiset $X_1$ is purely a subset of $\Z_{2N}$. We will use the same notation introduced in Section 2.3. (The results in this section are valid for all odd prime power $q$.)
	
Recall that 
\[
W_\cQ:=\{i\, (\mod{q^2+q+1}):Q(\omega^i)=0\}=\{d_0,d_1,\ldots, d_{q}\},
\]
which satisfies $W_\cQ\equiv 2(E_1\cup E_2)\,(\mod{q^2+q+1})$. 
For each $u \in W_{\mathcal Q}$ with $u\,(\mod{N})\in 2^{-1}I_1$, since 
	$W_{\mathcal Q}\equiv 2^{-1}(I_1\cup I_1 \cup I_2)\,(\mod{N})$, 
	there is exactly one $\ell_u\in \{1,2,\ldots,M-1\}$ such that $u+\ell_uN\in W_{\mathcal Q}$. Define 
	\[
	g_M(\omega^u)=\Tr_{q^3/q}(\omega^{2u+\ell_u N})\omega^{\ell_uN}. 
	\]
	\begin{lemma}\label{lem:nece:1}
		Let $\eta$ be the quadratic character of $\F_{q^3}$. Then, 
		$X_1$ is purely a subset in $\Z_{2N}$ if and only if $\eta(2)\not=\eta(g_M(\omega^u))$ for all $u\in W_{\mathcal Q}$ with 
		$u\,(\mod{N})\in 2^{-1}I_1$. 
	\end{lemma}
	\pro
	In the definition of ${\mathcal X}_{\mathcal Q}$ (see \eqref{eqn_defX00}), we take $d_0$ to be $u$. Then, by Lemma~\ref{lem:inva},  we have $2\omega^u,g_M(\omega^u)\omega^u\in {\mathcal X}_{\mathcal Q}$ or $2\omega^{u+q^2+q+1},g_M(\omega^u)\omega^{u+q^2+q+1}\in {\mathcal X}_{\mathcal Q}$. In either case, $X_1$ is purely a subset 
	of $\Z_{2N}$ if and only if $\eta(2)\not=\eta(g_M(\omega^u))$ for all $u\in W_{\mathcal Q}$ with 
	$u\,(\mod{N})\in 2^{-1}I_1$. 
	\qed
	
	\begin{lemma}\label{lem:comp:g} {\em (\cite[Lemma~5.8]{M18})}
		For $u\in W_{\mathcal Q}$ with 
		$u\,(\mod{N})\in 2^{-1}I_1$, 
		it holds that 
		\begin{equation}\label{eq:comp:g}
		\eta(g_M(\omega^u))=
		\eta(-1)\eta(1-\omega^{\frac{\ell_u (q+1)(q^3-1)}{M}})
		\eta(1-\omega^{\frac{2\ell_u q(q^3-1)}{M}}). 
		\end{equation}
	\end{lemma}
	
	In \cite{M18}, by using Lemmas~\ref{lem:nece:1} and \ref{lem:comp:g}, conditions on $q$ that guarantee $X_1$ is purely a subset of $\Z_{2N}$ were determined in the cases where $M=3$, $7$. We now generalize that result by using the following proposition. 
	
	\begin{proposition}\label{prop:comp:g} 
		Let $\eta$ be the quadratic character of $\F_{q^3}$. Then, 
		$X_1$ is purely a subset of $\Z_{2N}$ if and only if 
		$\eta(2)\ne \eta(1+\omega^\frac{\ell(q^3-1)}{M})$ for all $\ell\in \{1,2,\ldots,M-1\}$. 
	\end{proposition}
	\pro 
	First we note that
	\[
	\eta(1-\omega^{\frac{\ell_u(q+1)(q^3-1)}{M}})=\eta(1-\omega^{\frac{\ell_u(q^2+q)(q^3-1)}{M}})=\eta(1-\omega^{-\frac{\ell_u(q^3-1)}{M}})
	=\eta(-1)\eta(1-\omega^{\frac{\ell_u(q^3-1)}{M}}). 
	\]
	Simplifiying \eqref{eq:comp:g}, we obtain 
	\[
	\eta(g_M(\omega^u))=
	\eta(1+\omega^{\frac{\ell_u(q^3-1)}{M}}). 
	\]
	Let $L$ be the set of all $\ell_u$'s, i.e., $L=\{\ell_u\,:\,u,u+\ell_uN\in W_{\mathcal Q},u\,(\mod{N})\in 2^{-1}I_1\}$. Then, $|L|\le M-1$ since $|I_1|=\frac{M-1}{2}$. 
	
	We now show that $|L|=M-1$. 
	If $\ell_{u_1}=\ell_{u_2}$ for some $u_1,u_2\in W_{\mathcal Q}$, we have 
	$u_1,u_1+\ell_{u_1} N,u_2,u_2+\ell_{u_2} N\in W_{\mathcal Q}$ and 
	$u_1-(u_1+\ell_{u_1} N)=u_2-(u_2+\ell_{u_2} N)$, a contradiction to \eqref{eq:SingerW}.  
	Hence $L=\{1,2,\ldots,M-1\}$,
	and the conclusion of the proposition follows from Lemma~\ref{lem:nece:1}. 
	\qed
	\vspace{0.3cm}
	
	By Lemma~\ref{lem:nece:1} and Proposition~\ref{prop:comp:g}, we need to 
	study the question of whether $\Psi_{M,h,\alpha,\beta}$ is an infinite set, where 
	$\Psi_{M,h,\alpha,\beta}$ is defined in \eqref{def:psi}.

	\subsection{Multiplicative relations between the  $\left(1+\omega^\frac{\ell(q^3-1)}{M}\right)$'s}
	In this subsection,  in view of Proposition~\ref{prop:comp:g},  
	we study the question of whether $\eta(1+\omega^\frac{\ell(q^3-1)}{M})$, $1\le \ell\le M-1$, are all equal.  To this end, we first give 
	multiplicative relations between $1+\omega^\frac{\ell(q^3-1)}{M}$, $1\le \ell\le M-1$.  
	Let $\epsilon_M$ denote either $\omega^{\frac{q^3-1}{M}}$ or $\zeta_M=e^{2\pi i/M}$. 

For a positive inetegr $n$ let $\zeta_n$ be a primitive (complex) $n^{\rm th}$ root of unity, and let $D^{(n)}$ be the multiplicative group generated by $1-\zeta_n^i$, $i\in\{1,2,\ldots ,n-1\}$, modulo roots of unity. 
This group and its subgroups, in particular the subgroup of cyclotomic units, have been studied in \cite{C1,En1,F82,GK89,KR,La,R66,TS20}.
In particular,  nontrivial relations, called {\it Ennola relations}, between $(1-\zeta_n^i)$'s in $D^{(n)}$ were studied in \cite{C1,En1,KR}. 
The following relations  in $D^{(n)}$ are said to be {\it trivial relations}. 
\begin{itemize}
\item[1.] $1-\zeta_n^i=-\zeta_n^i(1-\zeta_n^{-i})$.  
\item[2.] A relation derived from the polynomial identity $\prod_{i=0}^{p-1}(1-x\zeta_p^i)=1-x^p$ for any prime $p$. 
\end{itemize}
Ennola~\cite{En1} gave for $n=105$ a relation which is not obtained from the 
trivial relations. Furthermore, he proved that twice any relation is a consequence of the trivial relations. 
Algorithms for finding  Ennola relations were studied in \cite{C1,KR}. Note that if we find an Ennola relation $R=1$ for  odd $n$, 
we have a relation on $(1+\zeta_n^i)$'s by dividing $\sigma_2(R)$ by $R$, where  
$R$ is a product of some integral power of $(1-\zeta_n^i)$'s and a root of unity, and $\sigma_2\in \Gal(\Q(\zeta_n)/\Q)$ is defined by $\sigma_2:\zeta_n\mapsto \zeta_n^2$. However, 
it is difficult to find Ennola relations in general. 
Furthermore, the signs or roots of unity involved in the multiplicative relations were usually ignored in those studies while here we need to be concerned about the signs. Moreover, our problem is weaker than that treated in the previous studies since we are interested in characterizing pairs $(q,M)$ such that   $\eta(1+\epsilon_M)=\cdots=\eta(1+\epsilon_M^{M-1})$.

	\begin{proposition}\label{cor:matome}
		Let $\eta$ be the quadratic character of $\F_{q^3}$ and $\ell$ be an integer, $1\le \ell \le M-1$ such that $\gcd{(\ell,M)}=1$. Then, the following hold: 
		\begin{itemize}
			\item[(1)] $\eta(1+\epsilon_M^{\ell})=\eta(1+\epsilon_M^{-\ell})$. 
			\item[(2)] 
			Let $s$ and $t$ be odd positive integers such that $t\not| s$. 
			Then, 
			$\eta(1+\epsilon_t^{s\ell})=
			\prod_{i=0}^{s-1}\eta(1+\epsilon_s^i\epsilon_t^\ell)$. In particular,  
			if $M=p^e$ is a prime power with $p$ a prime and $e\ge 2$,  $\prod_{i=0}^{p^j-1}\eta(1+\epsilon_{p^j}^i\epsilon_M^\ell)=
			\eta(1+\epsilon_M^{p^{j}\ell})$
			for $1\le j\le e-1$. 

			\item[(3)] $\prod_{i=0}^{\ord_M(2)-1}\eta(1+\epsilon_M^{2^i \ell})=1$. 
			\item[(4)] $\prod_{i=0}^{\ord_M(2)/2-1}\eta(1+\epsilon_M^{2^i\ell})=\eta(-1)$
			if $-1\in \langle 2\rangle\,(\mod{M})$. 
		\end{itemize}
	\end{proposition}
	\pro
	Note that $\eta(\epsilon_M^{\ell})=1$. 
The claims (1) and (2) follow from the trivial relations.  

We now prove the claim in (3). 
Since $2^{\ord_M(2)}\equiv 1\,(\mod{M})$, we have 
	\begin{align*}
	\prod_{i=0}^{\ord_M(2)-1}(1+\epsilon_M^{2^i \ell})&\,=\prod_{i=0}^{\ord_M(2)-1}
	\frac{1-\epsilon_M^{2^{i+1} \ell}}{1-\epsilon_M^{2^i \ell}}=
	\frac{1-\epsilon_M^{2^{\ord_M(2)} \ell}}{1-\epsilon_M^{\ell}}=1. 
	\end{align*}

Next, we prove the claim in (4). 
Since $2^{\ord_M(2)/2}\equiv -1\,(\mod{M})$, we have 
	\begin{align*}
	\prod_{i=0}^{\ord_M(2)/2-1}(1+\epsilon_M^{2^i \ell})&\,=\prod_{i=0}^{\ord_M(2)/2-1}
	\frac{1-\epsilon_M^{2^{i+1}\ell}}{1-\epsilon_M^{2^i \ell}}=
	\frac{1-\epsilon_M^{2^{\ord_M(2)/2}\ell}}{1-\epsilon_M^\ell}\\
	&\,=\frac{1-\epsilon_M^{-\ell}}{1-\epsilon_M^\ell}=
	-\epsilon_M^{-\ell} 
	\end{align*}
	The proof of the proposition is now complete. 
	\qed
	
	\begin{remark}\label{rem:others}
		In the case where $M$ is a prime power, say $M=p^e$, all 
{$\eta(1+\epsilon_M^{i})$ with $\gcd{(M,i)}\not=1$ are 
determined from $\eta(1+\epsilon_M^{i})$ with $\gcd{(M,i)}=1$ by 
Proposition~\ref{cor:matome}~(2).}  
In particular, if $\eta(1+\epsilon_M^{i})=\alpha$ for all $i$ with $\gcd{(i,M)}=1$, then $\eta(1+\epsilon_M^{i})=\alpha$ for all $i$ with $\gcd{(i,M)}\not=1$. 
	\end{remark}
	
	In the case where $M$ is not a prime power, we have the following nontrivial identity. 
	\begin{proposition}\label{prop:nece:mul}
		Let $M$ be a positive integer and $U_M=\{x\,:\,1\le x\le M-1,\gcd{(x,M)}=1\}$. 
		Assume that  $M$ is not a prime power, $-1\not \in \langle 2\rangle\,(\mod{M})$ and 
		$|\langle 2\rangle\,(\mod{M})|$ is even. 
		For any $\frac{\phi(M)}{4}$-subset $X$ of $U_M$ such that
		\begin{equation}\label{eq:x_x2x}
		X\cup -X\cup 2X \cup -2X =U_M
		\end{equation}
		and any integer $\ell$ with $\gcd{(\ell,M)}=1$, we have 
		\begin{equation}
		\prod_{x\in X}\eta(1+\epsilon_M^{2x \ell})=
		\eta(-1)^{\frac{\phi(M)}{4}+c}, 
		\end{equation}
		where 
		$c=|\{x\in X\cup 2X\,(\mod{M})\,:\,(M+1)/2\le x<M\}\}|$. 
	\end{proposition}
	\pro
	It is clear that 
	\begin{equation}\label{eq:m42}
	\prod_{x\in X}(1+\zeta_M^{2x \ell})=
	\prod_{x\in X}
	\frac{1-\zeta_M^{4x \ell}}{1-\zeta_M^{2x \ell}}
	=
	\frac{\prod_{x\in X\cup 2X}\zeta_{M}^{-x\ell}\prod_{x\in X\cup 2X}(1-\zeta_M^{2x \ell})}{\prod_{x\in X\cup 2X}\zeta_{M}^{-x\ell}\prod_{x\in X}(1-\zeta_M^{2x \ell})^2}. 
	\end{equation}
	We note that  $\prod_{x\in X\cup 2X}\zeta_{M}^{-x\ell}\prod_{x\in X\cup 2X}(1-\zeta_M^{2x \ell})=-1$ or $1$. Indeed, by noting that $\Phi_M(1)=1$ for 
	the $M$th cyclotomic polynomial $\Phi_M(x)$, we have 
	\begin{align*}
	\prod_{x\in X\cup 2X}\zeta_{M}^{-2x\ell}\prod_{x\in X\cup 2X}(1-\zeta_M^{2x \ell})^2=\prod_{x\in U_M}(1-\zeta_M^{2x \ell})=\Phi_M(1)=1. 
	\end{align*}
	On the other hand, 
	\begin{align*}
	\prod_{x\in X\cup 2X}\zeta_{M}^{-x\ell}\prod_{x\in X\cup 2X}(1-\zeta_M^{2x \ell})&\,=
	\prod_{x\in X\cup 2X}(\zeta_M^{x \ell}-\zeta_M^{-x \ell})\\
	&\,=(-4)^{|X|}\prod_{x\in X\cup 2X}\sin(2\pi x\ell/M)
	=(-4)^{\frac{\phi(M)}{4}}(-1)^{c_{X\cup 2X,\ell}} r
	\end{align*}
	for some positive $r\in \R$, where $c_{X\cup 2X,\ell}=|\{x\in X\cup 2X\,(\mod{M})\,:\,(M+1)/2\le x\ell<M\}\}|$. 
	Hence, we obtain $\prod_{x\in X\cup 2X}\zeta_{M}^{-x\ell}\prod_{x\in X\cup 2X}(1-\zeta_M^{2x \ell})=(-1)^{\frac{\phi(M)}{4}+c_{X\cup 2X,\ell}}$ and
	\begin{equation}\label{eq:444}
	\prod_{x\in X}(1+\zeta_M^{2x \ell})=
	\frac{(-1)^{\frac{\phi(M)}{4}+c_{X\cup 2X,\ell}}}{\prod_{x\in X\cup 2X}\zeta_{M}^{-x\ell}\prod_{x\in X}(1-\zeta_M^{2x \ell})^2}. 
	\end{equation}
	Next, we show that the parity of $c_{X\cup 2X,\ell}$ does not depend on the choice of $\ell$ (and also the choice of $X$ as a stronger claim). 
	Let $X'$ be an arbitrary union of cosets of $\langle 4\rangle\,(\mod{M})$ satisfying the condition 
	~\eqref{eq:x_x2x}. This is well defined since $-1\not\in \langle 2\rangle\,(\mod{M})$ and $|\langle 2\rangle\,(\mod{M})|$ is even. 
	Noting that $c_{X\cup 2X,\ell}=c_{\ell X \cup 2\ell X,1}$, it is enough to see that $c_{X\cup 2X,1}\equiv c_{X'\cup 2X',1}\,(\mod{2})$ for any $X$ satisfying the condition~\eqref{eq:x_x2x}. 
	
	Assume that $X'=\bigcup_{a\in A}a\langle 4\rangle$ for some $A\subseteq U_M$. For any $a\in A$, we have $a\langle 2\rangle\subseteq (X'\cup 2X')$ and 
	$-a\langle 2\rangle\cap (X'\cup 2X')=\emptyset$. On the other hand, we can write $-a\langle 2\rangle\cap (X\cup 2X)=\bigcup_{x\in S_a}\{x,2x\}$ for some $S_a\subseteq -a\langle 2\rangle$. By letting $Y=\bigcup_{a\in A}S_a$, we have $X'\cup 2X'=((X\cup 2X) \setminus (Y\cup 2Y)) \cup (-Y\cup -2Y)$. 
	Since $$c_{(X\cup 2X)\setminus \{x,2x\}\cup \{-x,-2x\},1}\in \{c_{X\cup 2X,1},c_{X\cup 2X,1}-2,c_{X\cup 2X,1}+2\}$$ for any $x\in X$, we have $$c_{X\cup 2X,1}\equiv c_{((X\cup 2X) \setminus (Y\cup 2Y)) \cup (-Y\cup -2Y),1}=c_{X'\cup 2X',1}\,(\mod{2}).$$ 
	
	Finally, since \eqref{eq:444} is also valid over $\F_{q}$, the conclusion of the proposition follows. 
	\qed
	\vspace{0.3cm}
	
	Let $\cP_M$ be the set of all prime powers $q$ such that $M\,|\,(q^2+q+1)$. 
	Define 
	\begin{align*}
	\Psi_{M,\alpha,\beta}&\,=\{q\in \cP_M\,:\,\eta(1+\epsilon_M^i)=\alpha, 1\le i<M,\eta(-1)=\beta\} 
	\end{align*}
	for $\alpha,\beta\in \{1,-1\}$. 
	In view of Proposition~\ref{prop:comp:g}, we consider when $\Psi_{M,\alpha,\beta}$ is empty. 
	\begin{proposition}\label{prop:-1in21}
		\begin{itemize}
			\item[(1)] If $-1\in \langle 2\rangle \,(\mod{M})$, then 
			$\Psi_{M,1,-1}=\emptyset$.
			\item[(2)] If $-1\in \langle 2\rangle \,(\mod{M})$ and $\ord_M(2)/2$ is even, then 
			$\Psi_{M,-1,-1}=\emptyset$.
			\item[(3)] If  $-1\in \langle 2\rangle \,(\mod{M})$ and $\ord_M(2)/2$ is odd, then 
			$\Psi_{M,-1,1}=\emptyset$.
			\item[(4)] If $\ord_M(2)$ is odd, then 
			$\Psi_{M,-1,\beta}=\emptyset$ for both $\beta=1,-1$.
		\end{itemize} 
	\end{proposition}
	\pro
	(1) If $q\equiv 3\,(\mod{4})$ and $-1\in \langle 2\rangle \,(\mod{M})$, by Proposition~\ref{cor:matome}~(4),  we have $\prod_{i=0}^{\ord_M(2)/2-1}\eta(1+\epsilon_M^{2^i})=-1$. 
	So it impossible to have $\eta(1+\epsilon_M^i)=1$ for all $i=1,2,\ldots,M-1$. 
	
	(2) If $q\equiv 3\,(\mod{4})$ and $-1\in \langle 2\rangle \,(\mod{M})$, by Proposition~\ref{cor:matome}~(4),  we have $$\prod_{i=0}^{\ord_M(2)/2-1}\eta(1+\epsilon_M^{2^i})=-1.$$ 
	On the other hand, if $\ord_M(2)/2$ is even and $\eta(1+\epsilon_M^i)=-1$ for all $i=1,2,\ldots,M-1$, we have $\prod_{i=0}^{\ord_M(2)/2-1}\eta(1+\epsilon_M^{2^i})=1$, a contradiction. 
	
	(3) If $q\equiv 1\,(\mod{4})$ and $-1\in \langle 2\rangle \,(\mod{M})$, by Proposition~\ref{cor:matome}~(4),   we have $$\prod_{i=0}^{\ord_M(2)/2-1}\eta(1+\epsilon_M^{2^i})=1.$$
	On the other hand, if $\ord_M(2)/2$ is odd and $\eta(1+\epsilon_M^i)=-1$ for all $i=1,2,\ldots,M-1$, we have $\prod_{i=0}^{\ord_M(2)/2-1}\eta(1+\epsilon_M^{2^i})=-1$, a contradiction. 
	
	(4)  By Proposition~\ref{cor:matome}~(3),  we have $\prod_{i=0}^{\ord_M(2)-1}\eta(1+\epsilon_M^{2^i})=1$. 
	On the other hand, if  $\eta(1+\epsilon_M^i)=-1$ for all $i=1,2,\ldots,M-1$, we have $\prod_{i=0}^{\ord_M(2)-1}\eta(1+\epsilon_M^{2^i})=-1$, a contradiction. 
	\qed
	\subsection{Determination of $\Psi_{M,h,\alpha,\beta}$ in the cases where $M=3,7,21$}
	%
	

	Below we determine prime powers $q$ such that $\eta(2)\not=\eta(1+\epsilon_M^\ell)$ for all $1\le \ell<M$ in the cases where $M=3,7,21$. Here the result in the case where $M=21$ is new. 
	\begin{proposition}{\em (\cite[Proposition~5.9]{M18})}\label{prop:33}
		Let $q\equiv 1\,(\mod{3})$ and $\epsilon_3\in \F_{q^3}$. Then, $\eta(2)\not=\eta(1+\epsilon_3^i)$ for $i=1,2$ if and only 
		if $q\equiv 7,13\,(\mod{24})$. 
	\end{proposition}
	\pro
	Since $1+\epsilon_3^2=-\epsilon_3$ and $1+\epsilon_3=-\epsilon_3^2$, we have $\eta(1+\epsilon_3)=\eta(1+\epsilon_3^2)=\eta(-1)=(-1)^{\frac{q-1}{2}}$. On the other hand, by the supplementary law of the quadratic reciprocity, we have 
	$\eta(2)=(-1)^{\frac{p^2-1}{8}}$. Hence, 
	$\eta(2)\not=\eta(1+\epsilon_3^i)$ for $i=1,2$ if and only 
	if $q\equiv 5,7\,(\mod{8})$. Combining with the assumption that $q\equiv 1\,(\mod{3})$, we see that the conclusion of the proposition follows. 
	\qed
	\begin{proposition}{\em (\cite[Proposition~5.10]{M18})}\label{prop:77}
		Let $q\equiv 2$ or $4\,(\mod{7})$ and $\epsilon_7\in \F_{q^3}$. Then, $\eta(2)\not=\eta(1+\epsilon_7^i)$ for $i=1,2,\ldots,6$ if and only 
		if $q\equiv 11,37,51,53\,(\mod{54})$. 
	\end{proposition}
	\pro
	Since $1+\epsilon_7^\ell=1+\epsilon_7^{q\ell}=1+\epsilon_7^{q^2\ell}$ and 
	$\prod_{i=0}^{2}\eta(1+\epsilon_7^{2^i \ell})=1$ for $\ell\in \{1,-1\}$ by Proposition~\ref{cor:matome}~(3), we have 
	$\eta(1+\epsilon_7^{i})=1$ for all $i=1,2,\ldots,6$. 
	Hence, 
	$\eta(2)\not=\eta(1+\epsilon_7^i)$ for $i=1,2,\ldots,6$ if and only 
	if $q\equiv 3,5\,(\mod{8})$ by the supplementary law of the quadratic reciprocity. 
	Combining with the assumption that $q\equiv 2,4\,(\mod{7})$, we see that the conclusion of the proposition follows. 
	\qed
	
	\begin{proposition}\label{prop:2121}
		Let $q\equiv 4$ or $16\,(\mod{21})$ and $\epsilon_{21}\in \F_{q^3}$. Then, $\eta(2)\not=\eta(1+\epsilon_{21}^i)$ for $i=1,2,\ldots,20$ if and only 
		if $q\equiv 37,109\,(\mod{168})$. 
	\end{proposition}
	\pro
	Let $X=\{1,4,16\}$. By Proposition~\ref{prop:nece:mul}, we have 
	$\prod_{i\in X}\eta(1+\epsilon_{21}^{i\ell})=\eta(-1)$ for any $\ell$ with $\gcd{(\ell,21)}=1$. On the other hand,
	$\eta(1+\epsilon_{21}^{\ell})=\eta(1+\epsilon_{21}^{4\ell})=
	\eta(1+\epsilon_{21}^{16\ell})$. 
	Hence $\eta(1+\epsilon_{21}^\ell)=\eta(-1)$ for all $\ell$ with 
	$\gcd{(\ell,21)}=1$. Furthermore, as in the proof of Proposition~\ref{prop:77}, we have $\eta(1+\epsilon_7^i)=1$ for  $i=1,2,\ldots,6$. Moreover, as in 
	the proof of Proposition~\ref{prop:33}, we have $\eta(1+\epsilon_3^i)=\eta(-1)$ for  $i=1,2$. 
	Hence, $\eta(2)\not=\eta(1+\epsilon_{21}^i)$ for $i=1,2,\ldots,20$ if and only 
	if $q\equiv 5\,(\mod{8})$ by the supplementary law of the quadratic reciprocity. Combining with the assumption that $q\equiv 4,16\,(\mod{21})$, we see that the conclusion of the proposition follows. 
	\qed
	\begin{remark}
		In the case where $M=21$, Proposition~\ref{cor:matome}~(1)--(4) imply the following relations between $\eta(1+\epsilon_{21}^i)$'s with $\gcd{(i,M)}=1$:  
		$\prod_{i\in \{1,2,4,8,16,11\}}\eta(1+\epsilon_{21}^{i\ell})=1$, 
		$\prod_{i\in \{1,8\}}\eta(1+\epsilon_{21}^{i\ell})=1$, 
		$\prod_{i\in \{1,4,10,13,16,19\}}\eta(1+\epsilon_{21}^{i\ell})=1$, 
		$\eta(1+\epsilon_{21}^{\ell})=\eta(1+\epsilon_{21}^{-\ell})$ 
		for any $\ell$ with $\gcd{(\ell,M)}=1$. Since it is impossible to derive $\prod_{i\in \{1,4,16\}}\eta(1+\epsilon_{21}^{i})=\eta(-1)$ from these relations,  Proposition~\ref{prop:nece:mul} is not a consequence of Proposition~\ref{cor:matome}~(1)--(4).  
	\end{remark}
	\section{Existence of primes in $\Psi_{M,h,\alpha,\beta}$: 
		Proof of Theorem~\ref{main:ext_thm1}}\label{sec:inf}
	Throughout this section, we assume that $M>3$ is an odd integer and there is an $h$, $1\le h\le M-1$, such that $M\,|\,(h^2+h+1)$.
      Let 
	$x_{i,j}=(-1)^j\sqrt{1+\zeta_M^i}$ for $i=0,1,\ldots,M-1$ and $j=0,1$, and  
	let 
	$Z_M=\{x_{i,j}\,:\,i=0,1,\ldots,M-1,j=0,1\}$. 
	Let $E_{M}$ denote the Galois extension of $\Q$ obtained by adjoining all elements of $Z_M$ and $\zeta_4$ to $\Q$. We now explain how to use Theorem~\ref{thm:chev} to prove Theorem~\ref{main:ext_thm1}. 

	Let ${\mathcal O}_{E_M}$ be the ring of integers 
	of $E_M$. For a fixed integer $1\le h\le M-1$ such that $M\,|\,(h^2+h+1)$, 
	let $p\equiv h\,(\mod{M})$ be an odd prime unramified in $E_M$ and $\Pe$ be a 
	prime ideal lying over $(p)$ in $E_M$. 
	Let $\pi_p$ be the Frobenius automorphism of 
	$({\mathcal O}_{E_M}/\Pe)/(\Z/(p))$, i.e., 
	$\pi_p:x\mapsto x^p$, and let $\sigma_\Pe$ be the Frobenius substitution with respect to $\Pe$ in $E_M/\Q$. Here, we note that ${\mathcal O}_{E_M}/\Pe$ contains the subfield $\{x+\Pe:x\in \Z[\zeta_M]\}$ of order $p^3$, 
	which is isomorphic to  $\Z[\zeta_M]/\pe$ for $\pe=\Pe\cap \Z[\zeta_M]$. Let $\eta$ be the quadratic character of 
	$\Z[\zeta_M]/\pe$. By using the correspondence between $\pi_p$ and $\sigma_{\Pe}$ and the facts that $\sqrt{1+\zeta_M^i}$ and $\zeta_4$ are units in $E_M$ and {$\Pe$ is lying over $p(\not=2)$,} it is straightforward to see that
	
	\begin{align}
	& \mbox{$p\equiv h\,(\mod{M})$, $\eta(1+\zeta_M^i)=\alpha$, $1\le i\le M-1$, }\nonumber\\
	& \mbox{\hspace{3cm} $\eta(2)=-\alpha$ and $\eta(-1)=\beta$  in $\Z[\zeta_M]/\pe$}\nonumber\\
	\mbox{$\Leftrightarrow$ }&
	\mbox{$\pi_p(\zeta_M+\Pe)=\zeta_M^h+\Pe$, $\pi_p^3(x_{i,0}+\Pe)=\alpha x_{i,0}+\Pe$, $1\le i\le M-1$, }\nonumber\\
	& \mbox{\hspace{3cm} $\pi_p^3(\sqrt{2}+\Pe)=-\alpha \sqrt{2}+\Pe$   and $\pi_p^3(\zeta_4+\Pe)=\beta\zeta_4+\Pe$ }\nonumber\\
	\mbox{$\Leftrightarrow$ }& 
	\mbox{$\sigma_{\Pe}(\zeta_M)=\zeta_M^h$, $\sigma_{\Pe}^3(x_{i,0})=\alpha x_{i,0}$, $1\le i\le M-1$,}\nonumber\\
	& \mbox{\hspace{3cm} $\sigma_{\Pe}^3(\sqrt{2})=-\alpha\sqrt{2}$  and 
		$\sigma_{\Pe}^3(\zeta_4)=\beta \zeta_4$.}\label{eq:equi}
	\end{align}
{Here, the conditions above are not depending on the choice of ${\Pe}$ lying over $p$.} 

	By this correspondence and Theorem~\ref{thm:chev} with $F=E_M$, $E=\Q$, $G=\Gal(E_M/\Q)$ and $\sigma_{\Pe}$ above, $\Psi_{M,h,\alpha,\beta}$ is 
	an infinite set if there is $\sigma\in G$ satisfying that 
	$\sigma(\zeta_M)=\zeta_M^h$, $\sigma^3(x_{i,0})=\alpha x_{i,0}$, $1\le i\le M-1$,  $\sigma^3(\sqrt{2})=-\alpha\sqrt{2}$  and 
	$\sigma^3(\zeta_4)=\beta \zeta_4$. 
	Thus we need to study the structure of $\Gal(E_M/\Q)$.

	We first list a few  basic properties of  $\Gal(E_M/\Q(\zeta_M))$ from Galois theory. 
	\begin{fact}\label{fact:11}{\em 
			\begin{itemize}
				\item[(i)] 
				$\Gal(E_M/\Q(\zeta_M))$ is a normal subgroup of $\Gal(E_M/\Q)$ since 
				$\Q(\zeta_M)$ is a normal extension of $\Q$. 
				\item[(ii)] $\Gal(E_M/\Q(\zeta_M))$ is an abelian $2$-group since all the extensions 
				$\Q(\zeta_M,x_{i,0})/\Q(\zeta_M)$, $i=0,1,\ldots,M-1$,  and 
				$\Q(\zeta_M,\zeta_4)/\Q(\zeta_M)$ are of degree at most $2$. 
				In particular, $\Gal(E_M/\Q(\zeta_M))$ is an elementary abelian $2$-group since it is isomorphic to  a subgroup of $ \prod_{i=0}^{M-1}\Gal(\Q(\zeta_M,x_{i,0})/\Q(\zeta_M))\times \Gal(\Q(\zeta_M,\zeta_4)/\Q(\zeta_M))$. 
		\end{itemize}}
	\end{fact}
	
	We next use the following fundamental result in  Kummer's theory. 
	See, e.g., \cite[Chapter~XI]{S86}.    
	
	\begin{theorem}\label{thm:kummer}
		Let $K$ be a field of characteristic $0$ containing $\zeta_n$ and $K^\ast=K\setminus \{0\}$ be the multiplicative group of $K$. 
		Let $R$ be a subgroup of $K^\ast$ containing $(K^\ast)^n=\{a^n\,:\,a\in K^\ast\}$, and $K(\sqrt[n]{R})$ denote the field extension of $K$ obtained by adjoining all elements in $\{\sqrt[n]{a}\,:\,a\in R\}$. Assume that $K(\sqrt[n]{R})/K$ is an abelian extension and $R/(K^\ast)^n$ is finite. Then $\Gal(K(\sqrt[n]{R})/K)$ is isomorphic to $R/(K^\ast)^n$. 
	\end{theorem}
	
	\begin{lemma}\label{lem:aa0}
Let $R=\langle 2,-1,1+\zeta_M^i,y\,:\,i=1,2,\ldots,M-1,y\in (\Q(\zeta_M)^\ast)^2\rangle$. Then, $\Gal(E_M/\Q(\zeta_M))$ is isomorphic to $R/(\Q(\zeta_M)^\ast)^2$. 
	\end{lemma}
	\pro
	Noting that $\Gal(E_M/\Q(\zeta_M))$ is abelian, apply Theorem~\ref{thm:kummer} with $K=\Q(\zeta_M)$ and  $n=2$. 
	\qed
	
	\begin{lemma}\label{lem:aa1}
		It holds that $2\not\in \langle -1,1+\zeta_M^i,y\,:\,i=1,2,\ldots,M-1,y\in (\Q(\zeta_M)^\ast)^2\rangle$.  
	\end{lemma}
	\pro 
	Assume that $2=(-1)^{i_0}\prod_{i=1}^{M-1}(1+\zeta_M^i)^{c_i}\cdot x^2$ for some $x\in \Q(\zeta_M)$. Since $(-1)^{i_0}\prod_{i=1}^{M-1}(1+\zeta_M^i)^{c_i}$ is a unit, {$x\in \Z[\zeta_M]$ and} $x$ is not a unit in $\Q(\zeta_M)$. 
	Thus $(2)$ is ramified in $\Q(\zeta_M)$. On the other hand, it is well known that $(2)$ is ramified in $\Q(\zeta_M)$ if and only if $4\,|\,M$, a contradiction.   
	\qed
	
	\begin{lemma}\label{lem:aa2}
		$-1\in \langle 1+\zeta_M^i,y\,:\,i=1,2,\ldots,M-1,y\in (\Q(\zeta_M)^\ast)^2\rangle$ if and only if $-1\in \langle 2\rangle \,(\mod{M'})$ for some divisor $M'>1$ of $M$.  
	\end{lemma}
	\pro 
	Assume that  $-1\not\in \langle 2\rangle \,(\mod{M'})$ for all divisors $M'>1$ of $M$, and $-1=\prod_{i=1}^{M-1}(1+\zeta_M^i)^{c_i}\cdot x^2$ for some integers $c_i$ and $x\in \Q(\zeta_M)$. 
	Let $\sigma\in \Gal(\Q(\zeta_M)/\Q)$ be defined by $\sigma:\zeta_M\mapsto \zeta_M^2$. Since $\ord_{M}(2)$ is odd, we have 
	\[
	\prod_{i=0}^{\ord_{M}(2)-1}\sigma^i(-1)=-1. 
	\]
	On the other hand, since $\prod_{i=0}^{\ord_{M}(2)-1}(1+\zeta_M^{2^i\ell})=1$ as in the proof of Proposition~\ref{cor:matome}~(3), 
	\[
	\prod_{i=0}^{\ord_{M}(2)-1}\Big(\prod_{j=1}^{M-1}\sigma^i(1+\zeta_M^j)^{c_j}\Big) \sigma^i(x)^2=\prod_{i=0}^{\ord_{M}(2)-1}\sigma^i(x)^2. 
	\]
	Hence, $-1$ is a square in $\Q(\zeta_M)$, which is impossible since $M$ is odd.  
	
	Conversely, if $-1\in \langle 2\rangle \,(\mod{M'})$ for some divisor $M'>1$ of $M$, we have $-1=\zeta_{M'}\prod_{i=0}^{\ord_{M'}(2)/2-1}\sigma^i(1+\zeta_{M'})$ 
 as in the proof of Proposition~\ref{cor:matome}~(4). 
	\qed
	\vspace{0.3cm}
	
{Note that the condition that $-1\not \in \langle 2\rangle \,(\mod{M'})$ for all divisors $M'>1$ of $M$ is equivalent to that the order of $2\in \Z$ is odd in  $(\Z/M\Z)^\times$. }
	Let $D_M$ be the field  obtained by adjoining 
	all $x_{i,0}$, $1\le i\le M-1$, to $\Q$. Furthermore, let $D_M'$ be the field obtained by adjoining $\zeta_4$ to $D_M$. 
 Summarizing the previous lemmas, we have the following.  
	\begin{proposition}\label{prop:matome2}
		\begin{itemize} 
\item[(1)] {$\Gal(E_M/\Q(\zeta_M))$ is isomorphic to 
			$\Gal(\Q(\zeta_M,\sqrt{2})/\Q(\zeta_M))\times 
			\Gal(D_M'/\Q(\zeta_M))$.} 
			\item[(2)] If {the order of $2\in \Z$ is odd in  $(\Z/M\Z)^\times$,}  
then $\Gal(E_M/\Q(\zeta_M))$ is isomorphic to 
			$\Gal(\Q(\zeta_M,\sqrt{2})/\Q(\zeta_M))\times \Gal(\Q(\zeta_M,\zeta_4)/\Q(\zeta_M))\times 
			\Gal(D_M/\Q(\zeta_M))$. 
		\end{itemize}
	\end{proposition}
{\pro Note that $\Gal(E_M/\Q(\zeta_M))$ is isomorphic to a subgroup of index at most $2$ of  $\Gal(\Q(\zeta_M,\sqrt{2})/\Q(\zeta_M))\times \Gal(D_M'/\Q(\zeta_M))$. By Lemmas~\ref{lem:aa0} and \ref{lem:aa1}, we have 
$$|\Gal(E_M/\Q(\zeta_M))|\not=|\Gal(D_M'/\Q(\zeta_M))|,$$ which implies that $\Gal(E_M/\Q(\zeta_M))$ is isomorphic to $\Gal(\Q(\zeta_M,\sqrt{2})/\Q(\zeta_M))\times \Gal(D_M'/\Q(\zeta_M))$.
Similarly, $\Gal(D_M'/\Q(\zeta_M))$ is isomorphic to a subgroup of index at most $2$ of  $\Gal(\Q(\zeta_M,\sqrt{-1})/\Q(\zeta_M))\times \Gal(D_M/\Q(\zeta_M))$. Then, by Lemmas~\ref{lem:aa0} and \ref{lem:aa2}, 
$\Gal(D_M'/\Q(\zeta_M))$ is isomorphic to $\Gal(\Q(\zeta_M,\zeta_4)/\Q(\zeta_M))\times 
			\Gal(D_M/\Q(\zeta_M))$ if the order of $2\in \Z$ is odd in  $(\Z/M\Z)^\times$. 
\qed}
\vspace{0.3cm}

Let $K$ be any subfield of  $E_M$.  
For any $\sigma\in \Gal(E_M/\Q)$, let $\sigma_{|K}$ denote the restriction of 
$\sigma$ to $K$. 
\begin{corollary}\label{cor:matome2}
		\begin{itemize} 
\item[(1)] {$\Gal(E_M/\Q)$ is isomorphic to 
			$\Gal(\Q(\sqrt{2})/\Q)\times 
			\Gal(D_M'/\Q)$. In particular, the group embedding $\varphi:\Gal(E_M/\Q)\to \Gal(\Q(\sqrt{2})/\Q)\times 
\Gal(D_M'/\Q)$ defined by $\varphi(\sigma)=(\sigma_{|\Q(\sqrt{2})},\sigma_{|D_{M}'})$ gives the isomorphism.} 
			\item[(2)] If the order of $2$ in $(\Z/M\Z)^\times$ is odd, then $\Gal(E_M/\Q)$ is isomorphic to  
			$\Gal(\Q(\sqrt{2})/\Q)\times \Gal(\Q(\zeta_4)/\Q)\times 
\Gal(D_M/\Q)$. {In particular, the group embedding $\varphi:\Gal(E_M/\Q)\to \Gal(\Q(\sqrt{2})/\Q)\times \Gal(\Q(\zeta_4)/\Q)\times \Gal(D_M/\Q)$ defined by $\varphi(\sigma)=(\sigma_{|\Q(\sqrt{2})},\sigma_{|\Q(\sqrt{-1})},
\sigma_{|D_{M}})$ gives the isomorphism.} 
\end{itemize}
	\end{corollary}
	\pro
	We only give a proof of Part (1) of the corollary.  Part (2) can be proved similarly. 
{As in the proof of Proposition~\ref{prop:matome2}, we have 
$E_M\not=D_M'$, which implies that $\Q(\sqrt{2})\cap D_M'=\Q$.  
Then, by Galois theory, $\varphi$ gives the isomorphism  between 
 $\Gal(E_M/\Q)$ and 
	$\Gal(\Q(\sqrt{2})/\Q)\times \Gal(D_M'/\Q)$.} 
	\qed
	\vspace{0.3cm}

	The following is our main theorem in this section. 
	
	\begin{theorem}\label{thm:main0}
		Let $h$ be a positive integer such that $1\le h\le M-1$ and $M\,|\,(h^2+h+1)$. Then the following hold: 
		\begin{itemize}
			\item[(1)]  $\Psi_{M,h,1,1}\cup \Psi_{M,h,1,-1}$ is an infinite set. 
			\item[(2)] If the order of $2$ in $(\Z/M\Z)^\times$ is odd, both $\Psi_{M,h,1,1}$ and $\Psi_{M,h,1,-1}$ are  infinite sets. 
		\end{itemize}
	\end{theorem}
	\pro 
	By \eqref{eq:equi} and Theorem~\ref{thm:chev}, we will show that there exists a 
	$\sigma\in \Gal(E_M/\Q)$ such that $\sigma(\zeta_M)=\zeta_M^h$,  $\sigma^3(x_{i,0})=x_{i,0}$, $1\le i\le M-1$, $\sigma^3(\sqrt{2})=-\sqrt{2}$ (and $\sigma^3(\zeta_4)=\beta \zeta_4$ for each $\beta\in \{1,-1\}$ for  the latter statement). 
	
	(1) Let $\sigma_1\in \Gal(\Q(\sqrt{2})/\Q)$ such that $\sigma_1(\sqrt{2})=-\sqrt{2}$. 
	Since 
	the exponent of $\Gal(D_M'/\Q(\zeta_M))$ is $2$, 
	there is a $\sigma_2\in \Gal(D_M'/\Q)$ such that  $\sigma_2(\zeta_M)=\zeta_M^h$ and $\sigma_2^3=\id$ or of order $2$. 
	If $\sigma_2^3$ is of order $2$, we use $\sigma_2^4$ in the place of $\sigma_2$. The newly named $\sigma_2$ satisfies 
	$\sigma_2(\zeta_M)=\zeta_M^h$ and $\sigma_2^3=\id$. 
	Therefore the element $\sigma=(\sigma_1,\sigma_2)\in \Gal(\Q(\sqrt{2})/\Q)\times \Gal(D_M'/\Q)\simeq \Gal(E_M/\Q)$ has the required property mentioned above.  
	
	(2) Let $\sigma_1$ be the same as in (1) above. 
	Furthermore, fix $\beta\in \{1,-1\}$ and let 
	$\sigma_2\in \Gal(\Q(\zeta_4)/\Q)$ be such that $\sigma_2(\zeta_4)=\beta\zeta_4$. 
	Similar to (1), 
	there is a $\sigma_3\in \Gal(D_M/\Q)$ such that $\sigma_3(\zeta_M)=\zeta_M^h$ and $\sigma_3^3=\id$. 
	Therefore the element
	$\sigma=(\sigma_1,\sigma_2,\sigma_3)\in \Gal(\Q(\sqrt{2})/\Q)\times \Gal(\Q(\zeta_4)/\Q)\times\Gal(D_M/\Q)\simeq \Gal(E_M/\Q)$ has the required property mentioned above. 
	\qed
	\vspace{0.3cm}
	
	{\bf Proof of Theorem~\ref{main:ext_thm1}: } In order to apply Theorem~\ref{thm:const_sum}, we require $q\equiv 3\pmod 4$; so we set $\beta=-1$, and consider $\Psi_{M,h,1,-1}$. By Theorem~\ref{thm:main0}~(2),  
	$\Psi_{M,h,1,-1}$ is an infinite set if the order of $2$ in $(\Z/M\Z)^\times$ is odd. Therefore by applying Theorem~\ref{thm:const_sum} 
	together with Theorem~\ref{thm:main_threevalued} and Proposition~\ref{prop:comp:g}, the conclusion of the theorem follows. \qed 
	
	\section{Concluding Remarks}\label{sec:conclu}
	In this paper, we gave a sufficient condition for the existence of strongly  regular Cayley graphs with parameters $(v,k,\lambda,\mu)=
	(p^{6},r(p^3+1),- p^3+r^2+3 r,r^2+r)$, where $r=(p^2-1)M/2$ and $M|(p^2+p+1)$. In particular, we proved  that strongly regular graphs with the above parameters exist for  sufficiently large primes $p\in \Psi_{M,h,\alpha,-1}$. Further, we considered the question for which $M,h$ and $\alpha$, $\Psi_{M,h,\alpha,-1}$ is an infinite set.  In fact, we proved that $\Psi_{M,h,\alpha,-1}$ is an infinite set if the order of $2$ in $(\Z/M\Z)^\times$ is odd. A few natural questions arise. 
	\begin{problem}
		Find  other classes of $M$ such that $\Psi_{M,h,\alpha,-1}$ is an infinite set.  
	\end{problem}
	\begin{problem}
		Evaluate the density of $\Psi_{M,h,\alpha,-1}$ when it is an infinite set. 
	\end{problem}
	We studied these problems in the case where $M$ is a prime power with certain extra conditions in the supplementary note~\cite{Mos}. However, the problems above are fully  open in the case where $M$ is not a prime power and the order of $2$ in $(\Z/M\Z)^\times$ is even except for the case $M=21$. 

\end{document}